\documentclass[graybox]{svmult}

\usepackage{type1cm}        % activate if the above 3 fonts are
\usepackage{makeidx}         % allows index generation
\usepackage{graphicx}        % standard LaTeX graphics tool
\usepackage{multicol}        % used for the two-column index
\usepackage[bottom]{footmisc}% places footnotes at page bottom

\usepackage{newtxtext}       % 
\usepackage{newtxmath}       % selects Times Roman as basic font

\usepackage{url}

\newenvironment{spmatrix}[1]
 {\def\mysubscript{#1}\mathop\bgroup\begin{pmatrix}}
 {\end{pmatrix}\egroup_{\textstyle\mathstrut\mysubscript}}

\usepackage{bm}
\usepackage{tikz}
\usepackage{wrapfig}
\usepackage[export]{adjustbox}
\usepackage{subcaption}
\usepackage[numbers]{natbib}
\usepackage{float}
\usepackage{xcolor}
\makeindex             % used for the subject index
\begin{document}

\title{Discontinuous Galerkin Model Order Reduction of Geometrically Parametrized Stokes Equation}
\author{Nirav Vasant Shah, Martin Wilfried Hess and Gianluigi Rozza}
\institute{Nirav Vasant Shah, \email{snirav@sissa.it} \at Scuola Internazionale Superiore di Studi Avanzati - via Bonomea, 265 - 34136 Trieste ITALY, 
\and Martin Wilfried Hess, \email{mhess@sissa.it} \at Scuola Internazionale Superiore di Studi Avanzati - via Bonomea, 265 - 34136 Trieste ITALY,
\and Gianluigi Rozza, \email{grozza@sissa.it} \at Scuola Internazionale Superiore di Studi Avanzati - via Bonomea, 265 - 34136 Trieste ITALY. }
\titlerunning{DG MOR of Geometrically Parametrized Stokes Equation}
\authorrunning{Shah et. al.}
\maketitle

\abstract*{The present work focuses on the geometric parametrization and the reduced order modeling of the Stokes equation. We discuss the concept of a parametrized geometry and its application within a reduced order modeling technique.  The full order model is based on the discontinuous Galerkin method with an interior penalty formulation. We introduce the broken Sobolev spaces as well as the weak formulation required for an affine parameter dependency. The operators are transformed from a fixed domain to a parameter dependent domain using the affine parameter dependency. The proper orthogonal decomposition is used to obtain the basis of functions of the reduced order model. By using the Galerkin projection the linear system is projected onto the reduced space. During this process, the offline-online decomposition is used to separate parameter dependent operations from parameter independent operations. Finally this technique is applied to an obstacle test problem.The numerical outcomes presented include experimental error analysis, eigenvalue decay and measurement of online simulation time.}

\abstract{The present work focuses on the geometric parametrization and the reduced order modeling of the Stokes equation. We discuss the concept of a parametrized geometry and its application within a reduced order modeling technique.  The full order model is based on the discontinuous Galerkin method with an interior penalty formulation. We introduce the broken Sobolev spaces as well as the weak formulation required for an affine parameter dependency. The operators are transformed from a fixed domain to a parameter dependent domain using the affine parameter dependency. The proper orthogonal decomposition is used to obtain the basis of functions of the reduced order model. By using the Galerkin projection the linear system is projected onto the reduced space. During this process, the offline-online decomposition is used to separate parameter dependent operations from parameter independent operations. Finally this technique is applied to an obstacle test problem.The numerical outcomes presented include experimental error analysis, eigenvalue decay and measurement of online simulation time.}

\section{Introduction}
\label{introduction}

Discontinuous Galerkin Method (DGM) has shown quite promising results for the elliptic problems~\cite{peraire} as well as for the hyperbolic problems~\cite{hyperbolic}. DGM uses polynomial approximation for sufficient accuracy and allows discontinuity at the interface for greater flexibility. Model Order Reduction (MOR) allows reducing the size of the system by retaining only ``dominant'' modes. The faster computations obtained by MOR has helped in many query contexts, real time computations and quick transfer of computational results to industrial problems. MOR in combination with geometric parametrization has emerged as an alternative to the shape optimization and has been used in many engineering applications.
As evident from above advantages, the application of geometric parametrization and reduced order modeling to discontinuous Galerkin method will remain at the forefront of scientific work. The present work is organized as follow. We first explain the concept of geometric parametrization. Thereafter, the governing equations, broken Sobolev spaces and weak formulation are stated. The affine expansion and Proper Orthogonal Decomposition (POD) are briefly described in the subsequent sections. Finally, an obstacle test problem demonstrates the application of the introduced method with outcomes involving comparison of full order and reduced order model solutions, error analysis and eigenvalue decay.

\section{Geometric parametrization}\label{geometric_parametrization_section}

Let us consider $\Omega = \Omega(\mu) \in \mathbb{R}^d$ as an open bounded domain. The parameter tuple $\mu \in \mathbb{P}$, where $\mathbb{P}$ is the parameter space, completely characterizes the domain. Also, consider a parameter tuple $\bar{\mu} \in \mathbb{P}$, as the known parameter tuple and $\Omega(\bar{\mu})$ as the reference domain, whose configuration is completely known. We divide the domain $\Omega(\mu)$ into $n_{su}$ triangular subdomains such that $\Omega(\mu) = \bigcup\limits_{i=1}^{n_{su}} \Omega_i(\mu) \ , \ \Omega_i(\mu) \bigcap \Omega_j(\mu) = \emptyset \ , \ \text{for} \ i \neq j$. The bijective mappings $\bm{F}_i(\cdot,\mu) : \Omega_i(\bar{\mu}) \rightarrow \Omega_i(\mu)$ link the reference subdomains $\Omega_i(\bar{\mu}) \subset \Omega(\bar{\mu})$ and the parametrized subdomains $\Omega_i(\mu) \subset \Omega(\mu)$. We consider here maps, $\bm{F}_i$, of the form,
\begin{gather*}
x = \bm{F}_i(\hat{x},\mu) = \bm{G}_{F,i}(\mu)\hat{x} + c_{F,i}(\mu) \ ; \\ \forall x \in \Omega_i(\mu) \ , \ \forall \hat{x} \in \Omega_i(\bar{\mu}) \ , \ \bm{G}_{F,i}(\mu) \in \mathbb{R}^{d \times d} \ , \ c_{F,i} \in \mathbb{R}^{d \times 1} \ , \ 1 \leq i \leq n_{su} \ .
\end{gather*}
The boundary of $\Omega(\mu)$, that is $\partial \Omega(\mu)$ is divided into a  Neumann boundary $\Gamma_N(\mu)$ and a Dirichlet boundary $\Gamma_D(\mu)$ i.e. $\partial \Omega(\mu) = \Gamma_N(\mu) \cup \Gamma_D(\mu)$. The Jacobian matrices $\bm{G}_{F,i}$ and the translational vectors $c_{F,i}$ depend only on parameter tuple $\mu$. The construction of maps $\lbrace \bm{F}_i \rbrace_{i=1}^{n_{su}}$ has been explained in literatures such as ~\cite{CRBM}.

\section{Discontinuous Galerkin formulation}
\label{DG_formulation}

The domain $\Omega$ is divided into $N_{el}$ number of triangular elements $\tau_k$ such that $\Omega = \bigcup\limits_{k=1}^{N_{el}} \tau_k$. The triangulation $\mathcal{T}$ is the set of all triangular elements i.e. $\mathcal{T} = \lbrace \tau_k \rbrace_{k=1}^{N_{el}}$. The internal boundary is denoted by $\Gamma = \bigcup\limits_{k=1}^{N_{el}} \partial \tau_k \backslash \partial \Omega$. $\overrightarrow{n}$ is the outward pointing normal to an edge of element.

The governing equations in strong form can be stated as,
\begin{flalign}\label{stokes_strong_form}
\begin{split}
\text{Stokes equation: } & -\nu \Delta \overrightarrow{u} + \nabla p = \overrightarrow{f} \ , \ \text{in } \Omega \ , \\
\text{Continuity equation: } & \nabla \cdot \overrightarrow{u} = 0 \ , \ \text{in} \ \Omega \ , \\
\text{Dirichlet condition: } & \overrightarrow{u} = \overrightarrow{u}_D \ , \ \text{on } \Gamma_D \ , \\
\text{Neumann condition: } & -p \overrightarrow{n} + \nu \overrightarrow{n} \cdot \nabla \overrightarrow{u} = \overrightarrow{t} \ , \ \text{on} \ \Gamma_N \ .
\end{split}
\end{flalign}

The velocity vector field $\overrightarrow{u}$ and pressure scalar field $p$ are the unknowns. $\nu$ is the material property known as kinematic viscosity. Vector $\overrightarrow{f}$ is the external force term or source term. $\overrightarrow{u}_D$ is the Dirichlet velocity and vector $\overrightarrow{t}$ is the Neumann value.

Let us introduce the broken Sobolev space, for any $p \in \mathbb{N}$,
\begin{equation*}
H^p(\Omega,\mathcal{T}) = \lbrace v \in L^2(\Omega) \ | \ v|_{\tau_k} \in H^p(\tau_k)  \ , \ \forall \tau_k \in \mathcal{T} \rbrace .
\end{equation*}

We consider finite dimensional subspaces of broken Sobolev spaces (see ~\cite{hyperbolic}), that is the spaces of discontinuous piecewise polynomial functions, for the unknowns.
\begin{equation*} \label{velocity_pressure_test}
\begin{split}
\text{For velocity: } \mathbb{V} = \lbrace \overrightarrow{\phi} \in (L^2(\Omega))^d | \ \overrightarrow{\phi} |_{\tau_k} \in (P^D(\tau_k))^d \ , \ \tau_k \in \mathcal{T} \rbrace \ , \\
\text{For pressure: } \mathbb{Q} = \lbrace \psi \in (L^2(\Omega)) | \ \psi |_{\tau_k} \in (P^{D-1}(\tau_k)) \ , \ \tau_k \in \mathcal{T} \rbrace \ .
\end{split}
\end{equation*}
Here, $P^D(\tau_k)$ denotes the space of polynomials of degree $D, \ D \geq 2$ over $\tau_k$. It is to be noted that, due to the application of interior penalty $(IP)$ and boundary penalty, the construction of subspace of Sobolev space is not required for imposing Dirichlet boundary condition.

In finite dimensional or discrete system, velocity approximation $\overrightarrow{u}_h(x)$ and pressure approximation $p_h(x)$ at any point $x \in \Omega$ are given by,
\begin{equation}\label{velocity_pressure_coefficients}
\overrightarrow{u}_h(x) = \sum\limits_{i=1}^{N_u} \overrightarrow{\phi}_i \hat{u}_i \ , \
p_h(x) = \sum\limits_{i=1}^{N_p} \psi_i \hat{p}_i \ ,
\end{equation}
where $\hat{u}_i$'s and $\hat{p}_i$'s are coefficients of velocity basis functions and pressure basis functions respectively. 

We expect that $\overrightarrow{u}_h \rightarrow \overrightarrow{u}$ and $p_h \rightarrow p$ as $N_u \rightarrow \infty$ and $N_p \rightarrow \infty$ respectively. Considering the scope of present work, the convergence analysis will not be discussed here. The readers are advised to refer to \cite{pacciarini},\cite{jump_mean_operator},\cite{riviere}.

In the subsequent sections, $\left( \cdot \right),\left( \cdot \right)_{\Gamma_D},\left( \cdot \right)_{\Gamma_N},\left( \cdot \right)_{\Gamma}$ represent the $L^2$ scalar product over $\Omega,\Gamma_D,\Gamma_N,\Gamma$ respectively. The jump operator $\left[ \cdot \right]$ and the average operator $\lbrace \cdot \rbrace$ are important concepts in the DGM formulation and are required to approximate the numerical flux. We use the jump and average operators as represented in \cite{jump_mean_operator}.

The weak form of the Stokes equation is given by,
\begin{gather}\label{stokes_weak_ch3}
a_{IP}(\overrightarrow{u},\overrightarrow{\phi}) + b(\overrightarrow{\phi},p) + \left( \lbrace p \rbrace,[\overrightarrow{n} \cdot \overrightarrow{\phi}] \right)_{\Gamma \cup \Gamma_D} = l_{IP}(\overrightarrow{\phi}) \ ,
\end{gather}
\begin{equation}
\begin{split}
a_{IP}(\overrightarrow{u},\overrightarrow{\phi}) = \left( \nabla \overrightarrow{u}, \nabla \overrightarrow{\phi} \right) + C_{11} \left( [\overrightarrow{u}],[\overrightarrow{\phi}] \right)_{\Gamma \cup \Gamma_D} \\ - \nu \left( \lbrace \nabla \overrightarrow{u}\rbrace ,[\overrightarrow{n} \otimes \overrightarrow{\phi}] \right)_{\Gamma \cup \Gamma_D} - \nu \left( [\overrightarrow{n} \otimes \overrightarrow{u}], \lbrace \nabla \overrightarrow{\phi} \rbrace \right)_{\Gamma \cup \Gamma_D} \ ,
\end{split}
\end{equation}
\begin{gather}
b(\overrightarrow{\phi},\psi) = -\int_{\Omega} \psi \nabla \cdot \overrightarrow{\phi} \ , \\
l_{IP}(\overrightarrow{\phi}) = \left( \overrightarrow{f},\overrightarrow{\phi} \right) + \left( \overrightarrow{t},\overrightarrow{\phi} \right)_{\Gamma_N} + C_{11} \left(\overrightarrow{u}_D,\overrightarrow{\phi}\right)_{\Gamma_D} - \left( \overrightarrow{n} \otimes \overrightarrow{u}_D, \nu \nabla \overrightarrow{\phi} \right)_{\Gamma_D} \ .
\end{gather}

The penalty parameter $C_{11}>0$ is an empirical constant to be kept large enough to maintain the coercivity of $a_{IP}(\overrightarrow{u},\overrightarrow{\phi})$ (see \cite{jump_mean_operator}).

The weak form of the continuity equation is as follows,
\begin{equation}\label{contiuity_weak_ch3}
\begin{split}
b(\overrightarrow{u},\psi) + ({\psi},[\overrightarrow{n} \cdot \overrightarrow{u}])_{\Gamma \cup \Gamma_D} = (\psi,\overrightarrow{n} \cdot \overrightarrow{u}_D)_{\Gamma_D} \ .
\end{split}
\end{equation}

In the discrete form the system of equations can be written as, 
\begin{equation} \label{Stokes_matrix_ch3}
\begin{spmatrix}{\textrm{Stiffness matrix}}
    \bm{A} & \bm{B} \\
    \bm{B}^T & \bm{0}
\end{spmatrix}
\begin{spmatrix}{\textrm{Solution vector}}
    U \\
    P
\end{spmatrix}
=
\begin{spmatrix}{\textrm{Right hand side (Known)}}
    F_1  \\
    F_2
\end{spmatrix}
\textrm{.}
\end{equation}

Here, $\bm{A}_{ij} = a_{IP} (\overrightarrow{\phi}_i,\overrightarrow{\phi}_j)$, $\bm{B}_{ij} = b(\overrightarrow{\phi}_i,\psi_j) + \left( \lbrace \psi_j \rbrace , [n \cdot \overrightarrow{\phi}_i]\right)_{\Gamma \cup \Gamma_D}$, $F_1 = l_{IP}(\overrightarrow{\phi}_i)$ and $F_2 = \left( \psi_j,\overrightarrow{n} \cdot \overrightarrow{u}_D \right)_{\Gamma_D}$ for $i=1,\ldots,N_u$ and $j=1,\ldots,N_p$. The column vectors $U$ and $P$ are coefficients $\hat{u}_i$'s and $\hat{p}_i$'s respectively (equation \eqref{velocity_pressure_coefficients}).

\section{Affine expansion}

We evaluate and solve the Stokes equation weak formulation on the reference domain $\Omega({\bar{\mu}})$. Given a parameter tuple $\mu \neq \bar{\mu}$, we need to evaluate the linear system of equations \eqref{Stokes_matrix_ch3} on a new domain $\Omega(\mu)$. To accomplish this, we use the affine expansion using linearity of equation and dividing $\Omega(\bar{\mu})$ into triangular subdomains $\Omega_i(\bar{\mu}) \ , \ i = \lbrace 1,2,\ldots,n_{su} \rbrace$ as explained earlier in the section \ref{geometric_parametrization_section}. The affine expansion of operators has been explained in the literatures such as \cite{CRBM}. The bilinear form $a_{IP}(\cdot,\cdot;\mu)$ can be expressed as,
\begin{gather}\label{affine_expansion_aip}
a_{IP}(\overrightarrow{u},\overrightarrow{\phi};\mu) = \sum\limits_{i=1}^{i=Q_a} \theta_a^i(\mu) a_{IP}^i(\overrightarrow{u},\overrightarrow{\phi};\bar{\mu}) \ ,
\end{gather}
for some finite $Q_a$ and some bilinear forms $\lbrace a_{IP}^i(\cdot,\cdot) \rbrace_{i=1}^{Q_a}$. The bilinear form $a_{IP}(\cdot,\cdot;\bar{\mu})$ is evaluated once on the reference domain $\Omega(\bar{\mu})$. To evaluate the bilinear form $a_{IP}(\cdot,\cdot;\mu)$ on the parametrized domain $\Omega(\mu)$, we use the affine expansion \eqref{affine_expansion_aip}. Since the evaluation of scalar terms $\lbrace \theta_a^i(\mu) \rbrace_{i=1}^{Q_a}$ is much faster than the evaluation of bilinear form $a_{IP}(\overrightarrow{u},\overrightarrow{\phi};\mu)$, significant speedup can be obtained with the help of affine expansion. Similar affine expansion can be used for other terms of the weak form \eqref{stokes_weak_ch3}. In the case of geometric parametrization, the affine expansion is essentially a change of variables \cite{geometric_para_2}. However, it is pertinent to explain two expansions as specific to DGM formulation.

\begin{itemize}
\item In order to transfer the terms containing jump and average operator, following approach is used in the present analysis.
\begin{equation*}\label{jump_average_term_split}
\begin{split}
\left(\lbrace \nabla \overrightarrow{\phi} \rbrace , \left[ \overrightarrow{n} \otimes \overrightarrow{\phi}  \right]  \right) = \left( \nabla \overrightarrow{\phi}^+ , \overrightarrow{n}^+ \otimes \overrightarrow{\phi}^+ \right) + \left( \nabla \overrightarrow{\phi}^+ , \overrightarrow{n}^- \otimes \overrightarrow{\phi}^- \right) + \\ 
\left( \nabla \overrightarrow{\phi}^- , \overrightarrow{n}^+ \otimes \overrightarrow{\phi}^+ \right) + \left( \nabla \overrightarrow{\phi}^- , \overrightarrow{n}^- \otimes \overrightarrow{\phi}^- \right) \ .
\end{split}
\end{equation*}
Each term on the right hand side of the above equation can be transformed using the affine map. 

\item The coercivity term $C_{11}\left( [\overrightarrow{\phi}],[\overrightarrow{u}] \right)_{\Gamma \cup \Gamma_D}$ is not transformed but used as evaluated on reference domain $\Omega(\bar{\mu})$. The affine transformation is given by,
\begin{equation*}
\begin{split}
C_{11}\left( [\overrightarrow{\phi}(x),\overrightarrow{u}(x)] \right)_{\Gamma(\mu) \cup \Gamma_D(\mu)} = C_{11} \alpha \left( [\overrightarrow{\phi}(\bm{F}(\hat{x})),\overrightarrow{u}(\bm{F}(\hat{x}))] \right)_{\Gamma(\bar{\mu}) \cup \Gamma_D(\bar{\mu})} \ , \\
\alpha = \frac{\text{length of }\left( \Gamma(\mu) \cup \Gamma_D(\mu)\right)}{\text{length of }\left( \Gamma(\bar{\mu}) \cup \Gamma_D(\bar{\mu})\right)} \ , \ \hat{x} \in \Omega(\bar{\mu}) \ , \ x \in \Omega(\mu) \ .
\end{split}
\end{equation*}
Since, $C_{11}$ is an empirical coefficient replacing $C_{11} \alpha$ with $C_{11}$ will not change the formulation as long as the coercivity of $a_{IP}(\overrightarrow{u},\overrightarrow{\phi}) $ over parameter space $\mathbb{P}$ is maintained.
\end{itemize}

\section{Reduced basis method}\label{rb_section}

Snapshot POD exploits the information contained in the snapshots to construct low dimensional reduced basis space which can approximate the solution within desirable accuracy. The offline phase consists of construction of reduced basis space while the online phase consists of computing coefficients of the reduced basis. For detailed explanation about POD-Galerkin method and offline-online decomposition, we refer to ~\cite{CRBM}.

As first step, the DGM solutions based on $\mu_n, n \in \lbrace 1,....,n_s \rbrace$ are calculated i.e. $n_s$ snapshots are generated. The velocity snapshots and the pressure snapshots are stored in $\bm{S}_v \in \mathbb{R}^{N_u \times n_s}$ and $\bm{S}_p \in \mathbb{R}^{N_p \times n_s}$ respectively. Let us also introduce inner product matrices $\bm{M}_v \in \mathbb{R}^{N_u \times N_u}$ and $\bm{M}_p \in \mathbb{R}^{N_p \times N_p}$.

\begin{gather*}
\bm{M}_{v,ij} = \int_{\Omega} \overrightarrow{\phi}_i \cdot \overrightarrow{\phi}_j + \sum_{k=1}^{N_{el}} \int_{\tau_k} \nabla \overrightarrow{\phi}_i : \nabla \overrightarrow{\phi}_j \ , \ i,j = 1, \ldots, N_u \ , \\
\bm{M}_{p,ij} = \int_{\Omega} \psi_i \psi_j \ , \ i,j = 1, \ldots, N_p \ .
\end{gather*}

The dimension of the reduced basis is denoted as $N$ and it is asserted that $N << N_u, \ N < n_s$. Proper Orthogonal Decomposition obtains orthogonal basis for the low dimensional reduced basis space, by using spectral decomposition.
\begin{equation}\label{snapshot_eigen_value}
\bm{S}_v^T \bm{M}_v \bm{S}_v = \bm{V} \bm{\Theta} \bm{V}^T \ .
\end{equation}
The columns of $\bm{V}$ are eigenvectors and $\Theta$ has eigenvalues $\theta_i \ , \ 1 \leq i,j \leq n_s$, in sorted order ($\theta_1 \geq \ldots \geq \theta_{n_s}$) such that, $\Theta_{ij} = \theta_i \delta_{ij}$. Eigenvalue decay, the drop in the magnitude of the eigenvalues, provides upper bound for the error between the solution computed by full order model and the solution computed by POD (see ~\cite{CRBM}).

The projection matrix $\bm{B}_v \in \mathbb{R}^{N_u \times N}$, used for the projection from the space of full order model to the space of reduced order model, is given by, 
\begin{equation}
\bm{B}_v = \bm{S}_v \bm{V} \bm{\Theta}^{-\frac{1}{2}} \bm{R} \ , \ \bm{R} = [\bm{I}_{N \times N} ; \bm{0}_{(n_s-N) \times N}] \ ,
\end{equation}
where, $\bm{I}_{N \times N}$ is the identity matrix of size $N \times N$.
The reduced basis space $\bm{B}_p$ can be generated in a similar manner using the pressure snapshots $\bm{S}_p$ and the inner product matrix $\bm{M}_p$. Above procedure is performed during the offline phase.

The discrete system of equations is projected onto the reduced basis space by Galerkin projection as,
\begin{equation} \label{Stokes_matrix_reduced}
\begin{spmatrix}{\tilde{K}}
    \bm{B}_v^T \bm{A}(\mu) \bm{B}_v & \bm{B}_v^T \bm{B}(\mu) \bm{B}_p \\
    \bm{B}_p^T \bm{B}(\mu)^T \bm{B}_v & \bm{0}
\end{spmatrix}
\begin{spmatrix}{\zeta}
    U_N \\
    P_N
\end{spmatrix}
=
\begin{spmatrix}{\tilde{F}}
    \bm{B}_v^T F_1(\mu)  \\
    \bm{B}_p^T F_2(\mu)
\end{spmatrix} \ .
\end{equation}
The solution vectors $U$ and $P$ (equation \eqref{Stokes_matrix_ch3}) are then computed as $U = \bm{B}_v U_N \ , \ P = \bm{B}_p P_N$. Projection onto the reduced basis space, solution of smaller system of equations and computation of $U$ and $P$ are steps performed during online phase. During the online phase, the matrices $\bm{A}(\mu)$, $\bm{B}(\mu)$ and the vectors $F_1(\mu)$, $F_2(\mu)$ are evaluated using affine expansion.

\section{A numerical example}

The numerical experiments were performed using RBmatlab ~\cite{rbmatlab},~\cite{master_thesis}. The reference domain $\Omega({\bar{\mu}})$ is the unit square domain $[0,1] \times [0,1]$ with triangle having vertices $(0.3,0),(0.5,0.3),(0.7,0)$ as obstacle. The domain $\Omega(\bar{\mu})$ is divided into $9$ mutually non-overlapping subdomains. Two geometric parameters, the coordinates of the tip of the obstacle, with reference values collected in parameter tuple $\bar{\mu} = (0.5,0.3)$ characterize the domain. The $x-$direction refers to the horizontal direction and the $y-$direction refers to the vertical direction. The boundary ${x=0}$ is a Dirichlet boundary with inflow velocity at point $(0,y)$ as $u = (y(1-y), 0)$. The boundary ${x = 1}$ is a Neumann boundary with zero Neumann value i.e. $\overrightarrow{t} = (0, 0)$. Other boundaries are Dirichlet boundary with no slip condition. The source term is $\overrightarrow{f} = (0,0)$.

The training set contained $100$ uniformly distributed random parameters within the $[0.4,0.6] \times [0.2,0.4]$. The test set contained $10$ uniformly distributed random parameters within the range $[0.4,0.6] \times [0.2,0.4]$. For velocity basis function polynomial of degree $D = 2$ and for pressure basis function polynomial of degree ${D-1} = 1$ were used. The number of velocity degrees of freedom and pressure degrees of freedom were $N_u = 4704$ and $N_p = 1176$ respectively.

\begin{figure}[H] %[t!] % "[t!]" placement specifier just for this example
\begin{subfigure}{0.31\textwidth}
\includegraphics[width=\linewidth]{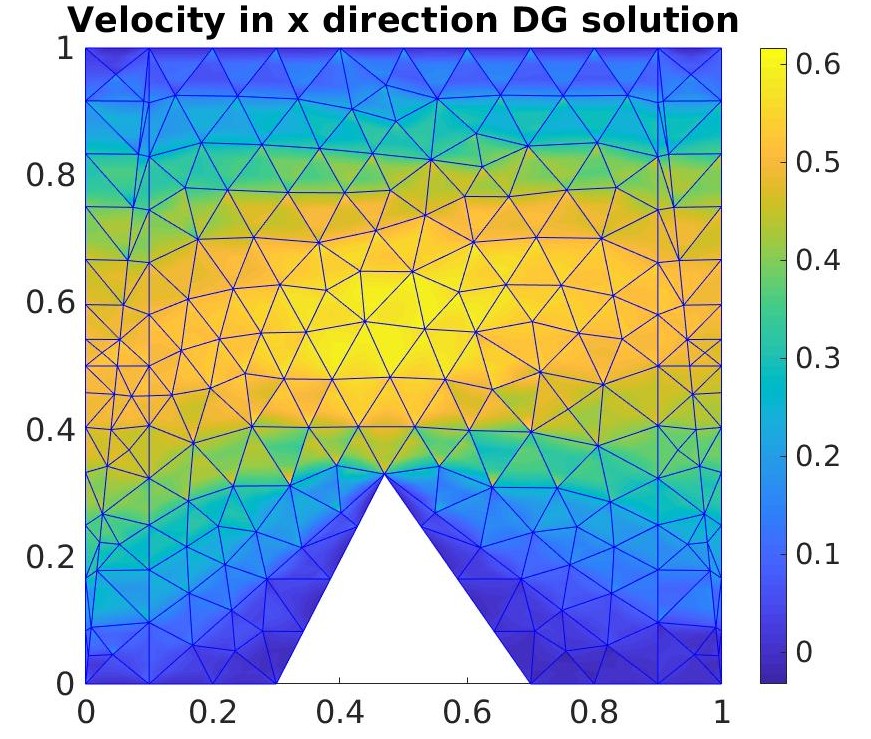}
\caption{Velocity $x-$direction DGM solution} \label{vel_x_dg}
\end{subfigure}\hspace*{\fill}
\begin{subfigure}{0.31\textwidth}
\includegraphics[width=\linewidth]{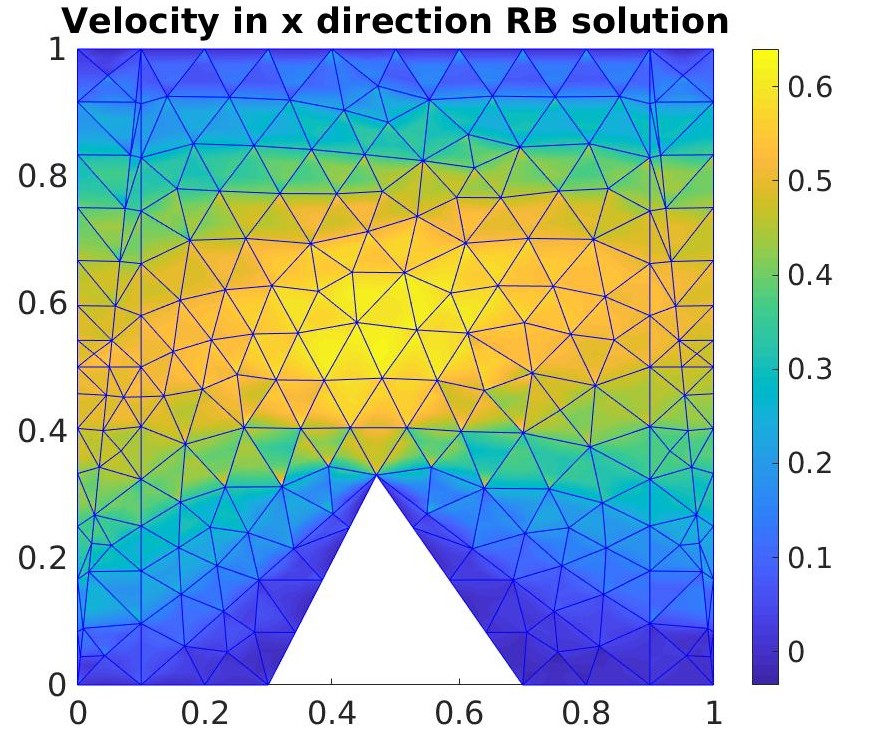}
\caption{Velocity $x-$direction RB solution} \label{vel_x_rb}
\end{subfigure}\hspace*{\fill}
\begin{subfigure}{0.31\textwidth}
\includegraphics[width=\linewidth]{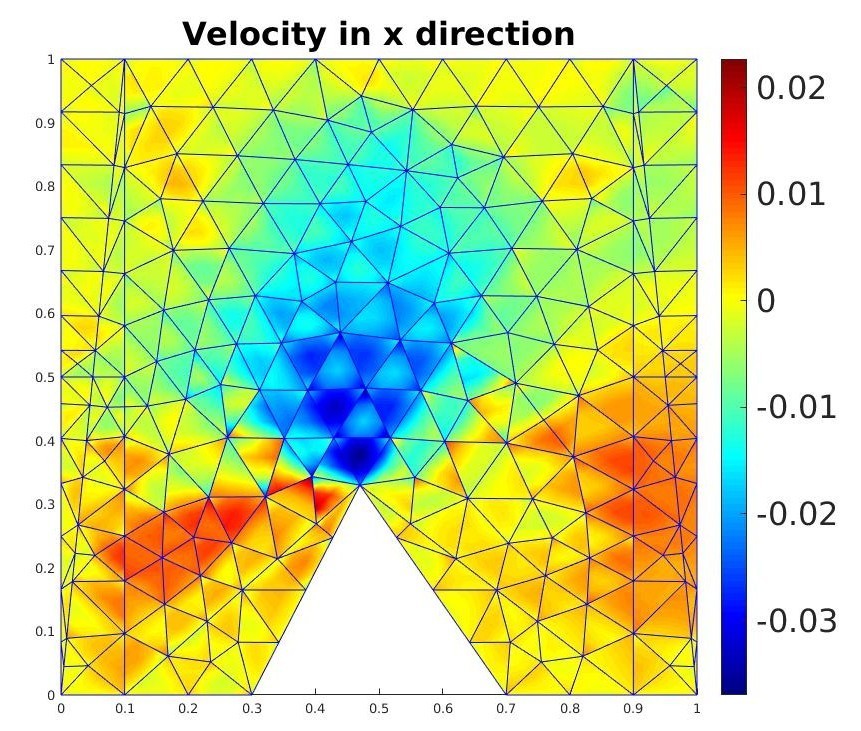}
\caption{$x-$component of Velocity absolute error $\overrightarrow{u}_h-\overrightarrow{u}_N$} \label{error_x_vel}
\end{subfigure}

\begin{subfigure}{0.31\textwidth}
\includegraphics[width=\linewidth]{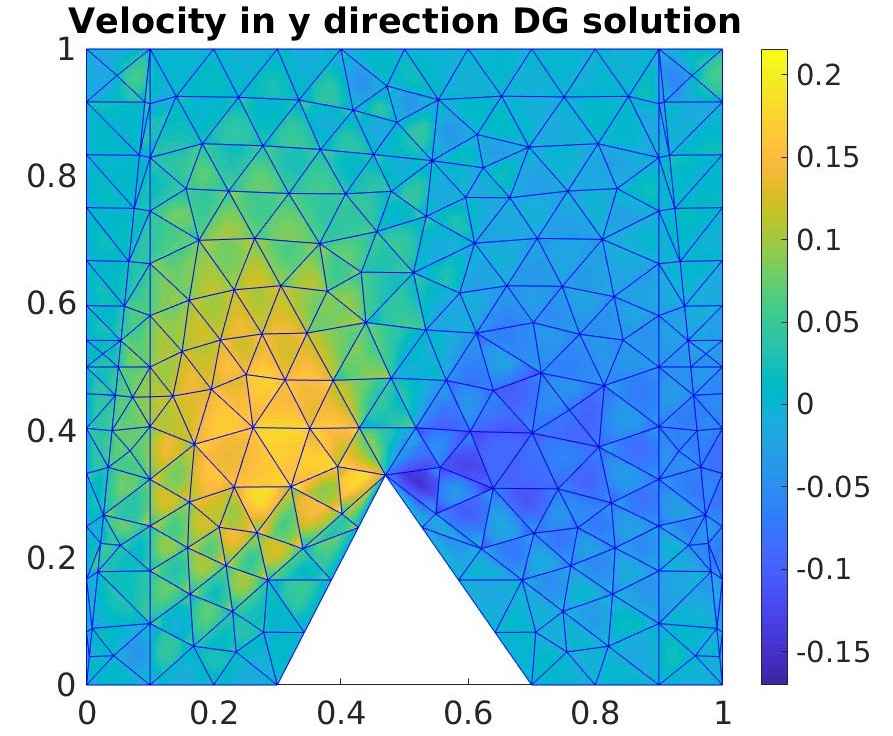}
\caption{Velocity $y-$direction DGM solution} \label{vel_y_dg}
\end{subfigure}\hspace*{\fill}
\begin{subfigure}{0.31\textwidth}
\includegraphics[width=\linewidth]{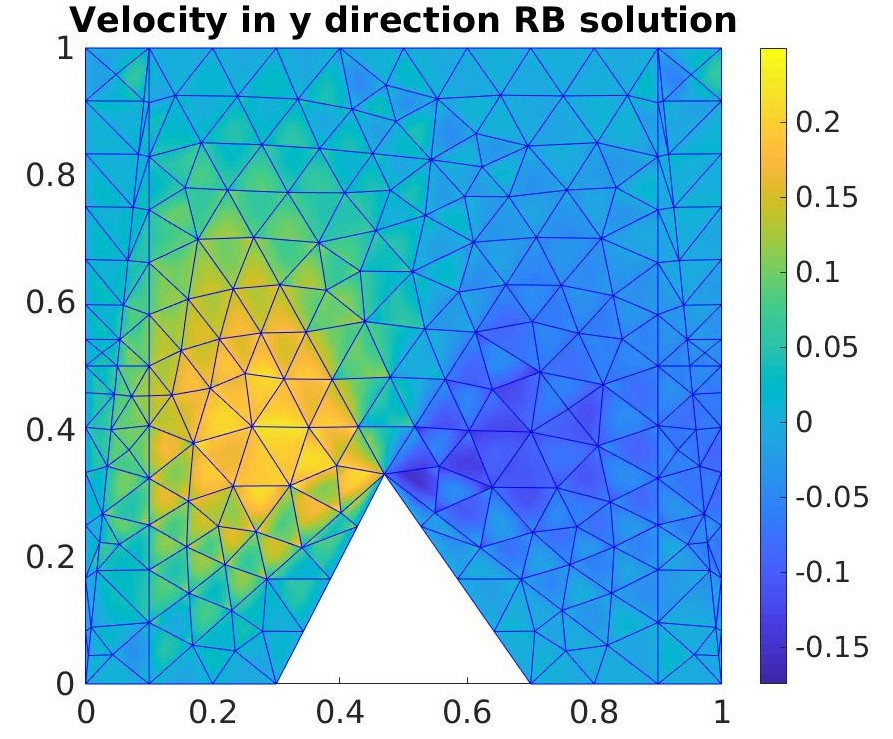}
\caption{Velocity $y-$direction RB solution} \label{vel_y_rb}
\end{subfigure}\hspace*{\fill}
\begin{subfigure}{0.31\textwidth}
\includegraphics[width=\linewidth]{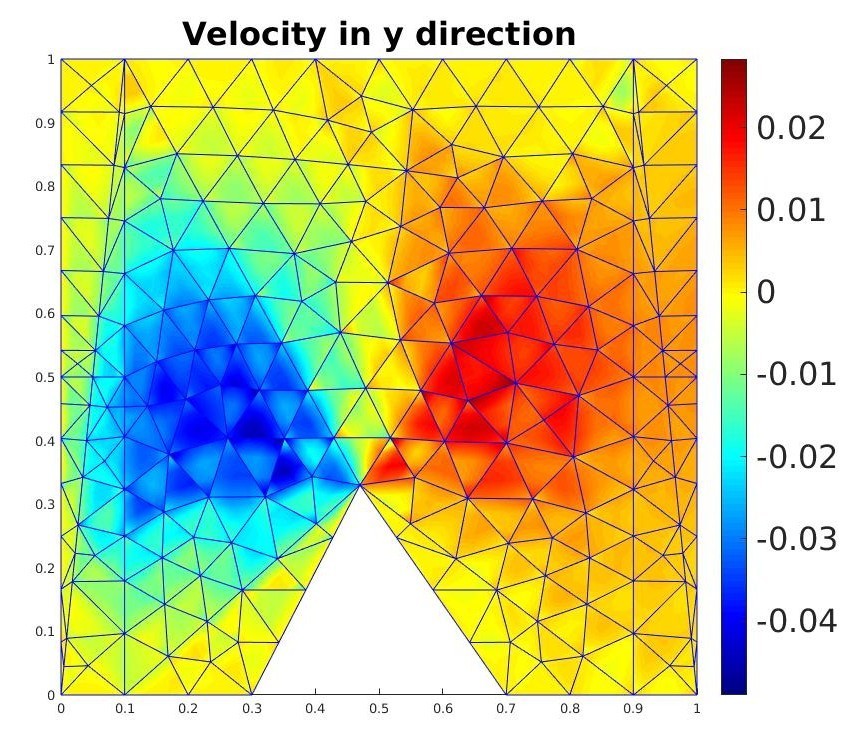}
\caption{$y-$component of Velocity absolute error $\overrightarrow{u}_h-\overrightarrow{u}_N$} \label{error_y_vel}
\end{subfigure}

%\begin{subfigure}{0.31\textwidth}
%\includegraphics[width=\linewidth]{offline_pressure_at_47_33.jpg}
%\caption{Pressure DGM solution} \label{pre_dg}
%\end{subfigure}\hspace*{\fill}
%\begin{subfigure}{0.31\textwidth}
%\includegraphics[width=\linewidth]{online_pressure_at_47_33.jpg}
%\caption{Pressure RB solution} \label{pre_rb}
%\end{subfigure}
%\begin{subfigure}{0.31\textwidth}
%\includegraphics[width=\linewidth]{pressure_error_at_47_33.jpg}
%\caption{Pressure absolute error $p_h-p_N$} \label{pre_error}
%\end{subfigure}
\caption{DGM and RB solution $\mu = (0.47,0.33)$} 
\label{dg_rb_solution_47_33}
\end{figure}

\begin{figure}[H]
\begin{subfigure}{0.48\textwidth}
\includegraphics[width=\linewidth]{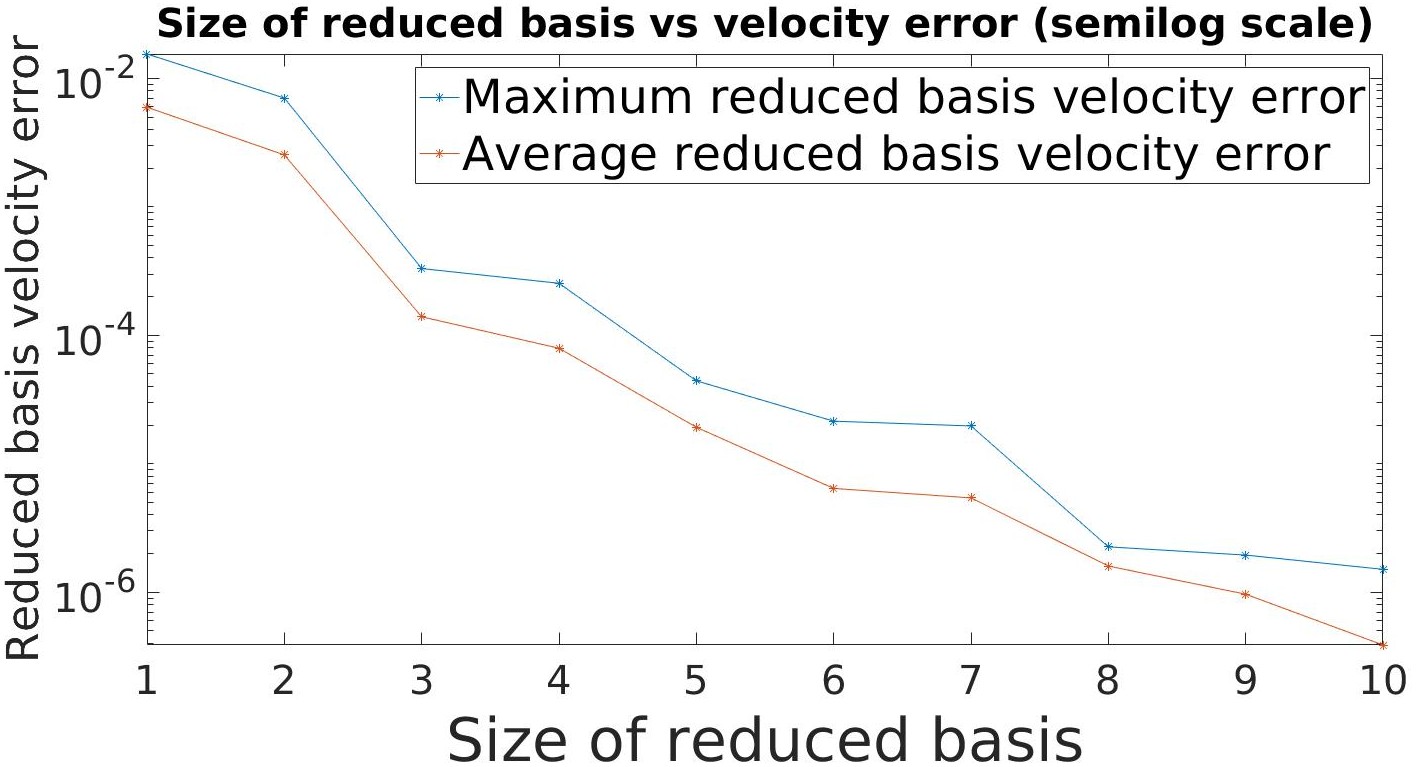}
\caption{Size of the reduced basis space vs. Relative error in velocity with inner product induced by $\bm{M}_v$} \label{error_vs_basis_velocity}
\end{subfigure}\hspace*{\fill}
\begin{subfigure}{0.48\textwidth}
\includegraphics[width=\linewidth]{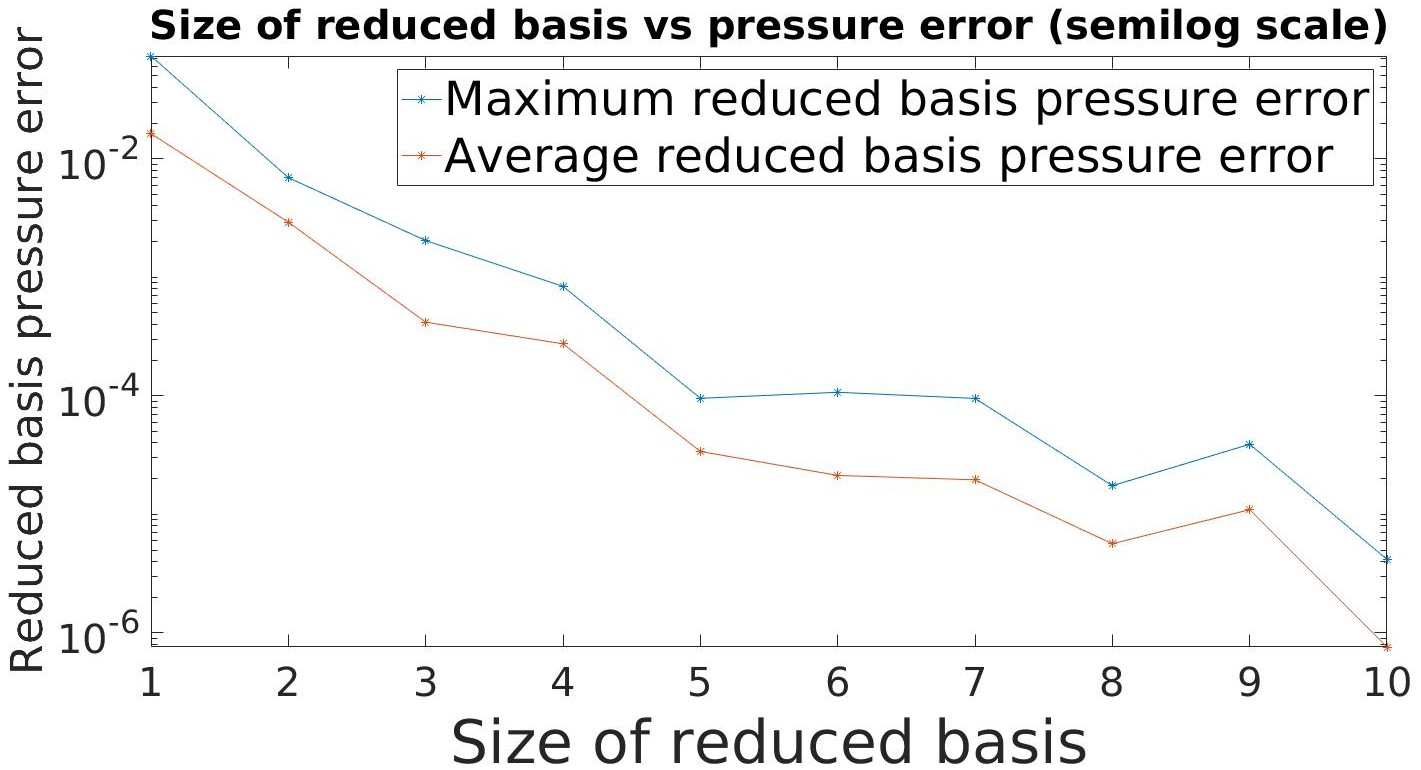}
\caption{Size of the reduced basis space vs. Relative error in pressure with inner product induced by $\bm{M}_p$} \label{error_vs_basis_pressure}
\end{subfigure}
  \caption{Size of the reduced basis space vs Relative error} 
\label{error_vs_basis}
\end{figure}

\begin{figure}[H]
\begin{subfigure}{0.31\textwidth}
\includegraphics[width=\linewidth]{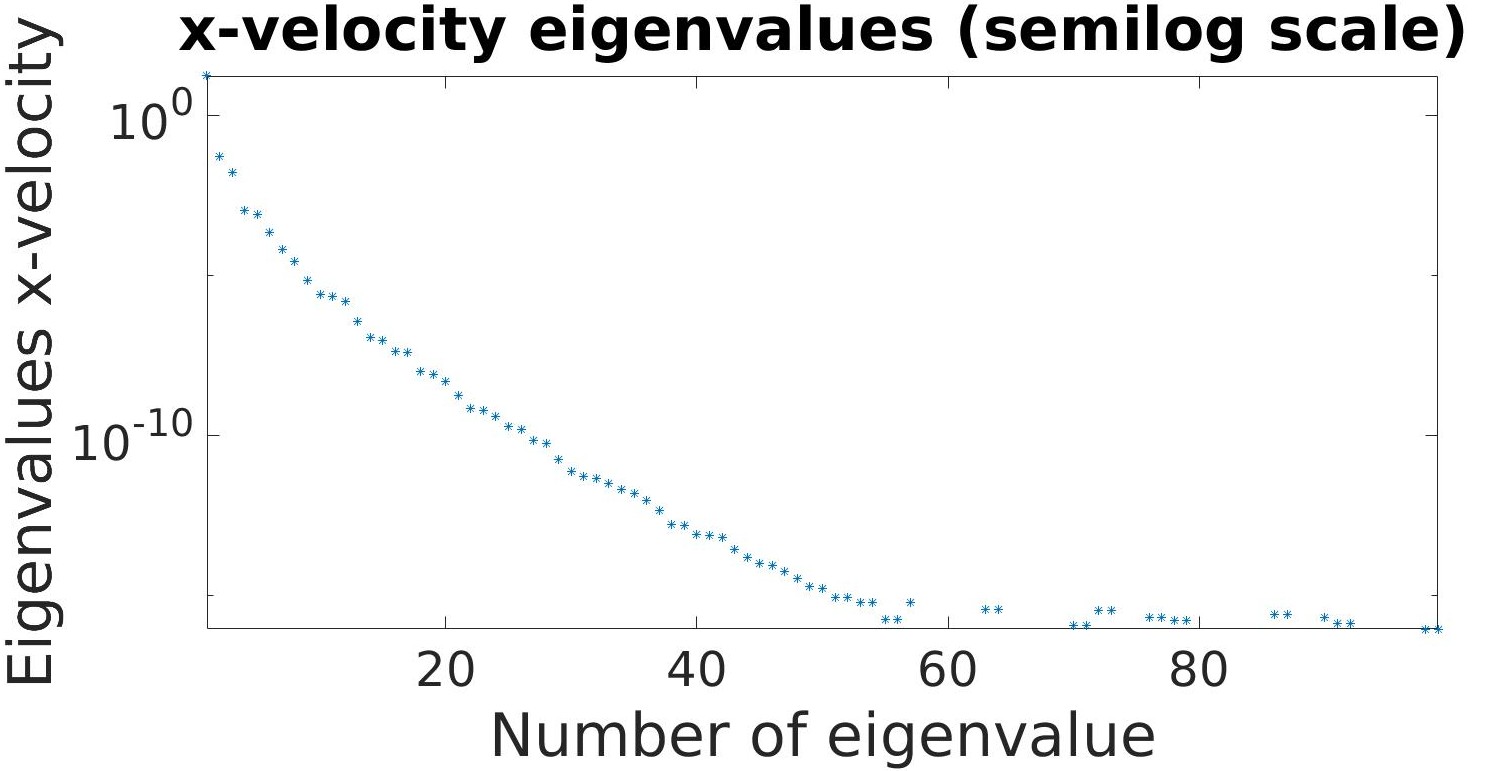}
\caption{$x-$Velocity eigenvalues (semilog scale)} \label{vel_x_ev}
\end{subfigure}\hspace*{\fill}
\begin{subfigure}{0.31\textwidth}
\includegraphics[width=\linewidth]{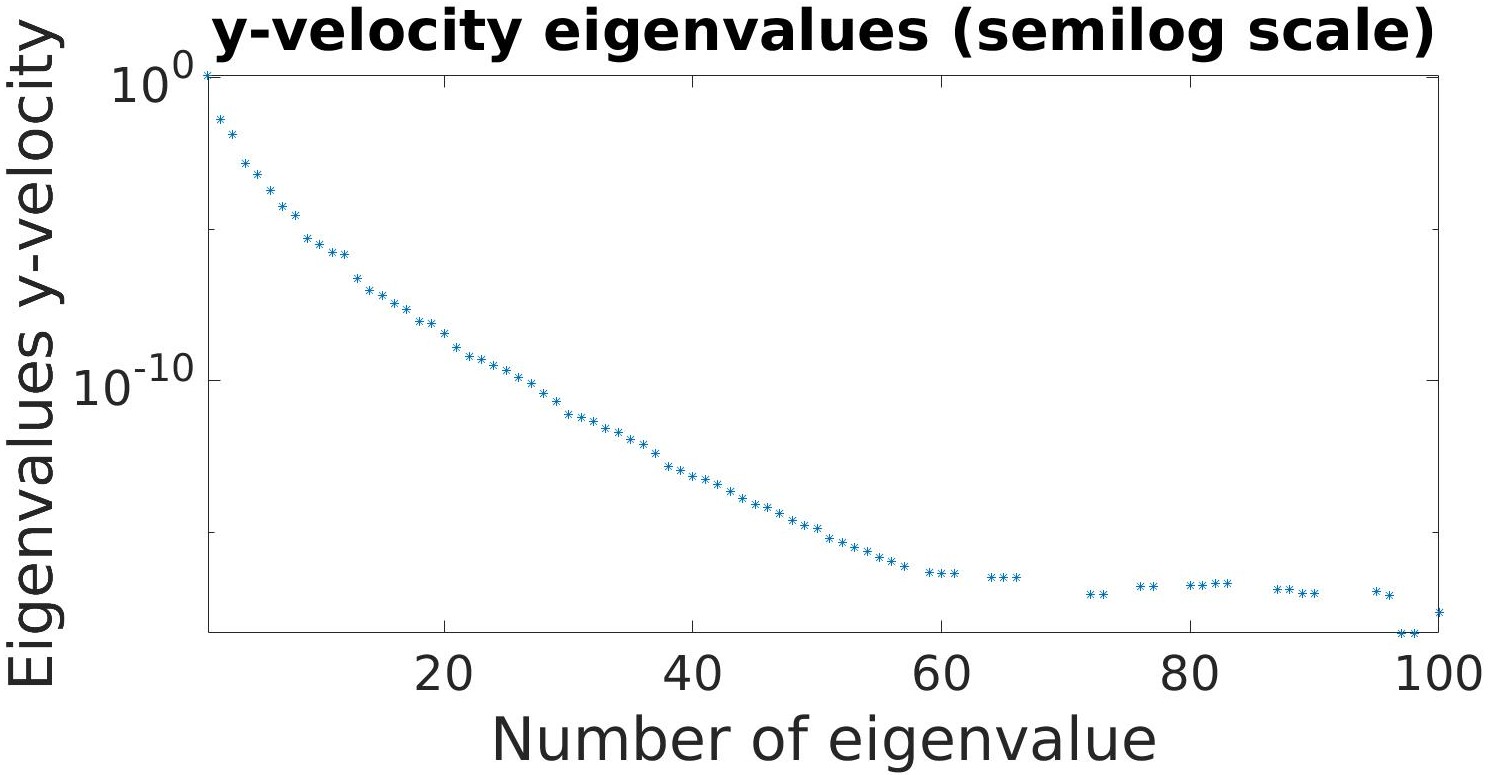}
\caption{$y-$Velocity eigenvalues (semilog scale)} \label{vel_y_ev}
\end{subfigure}
\begin{subfigure}{0.31\textwidth}
\includegraphics[width=\linewidth]{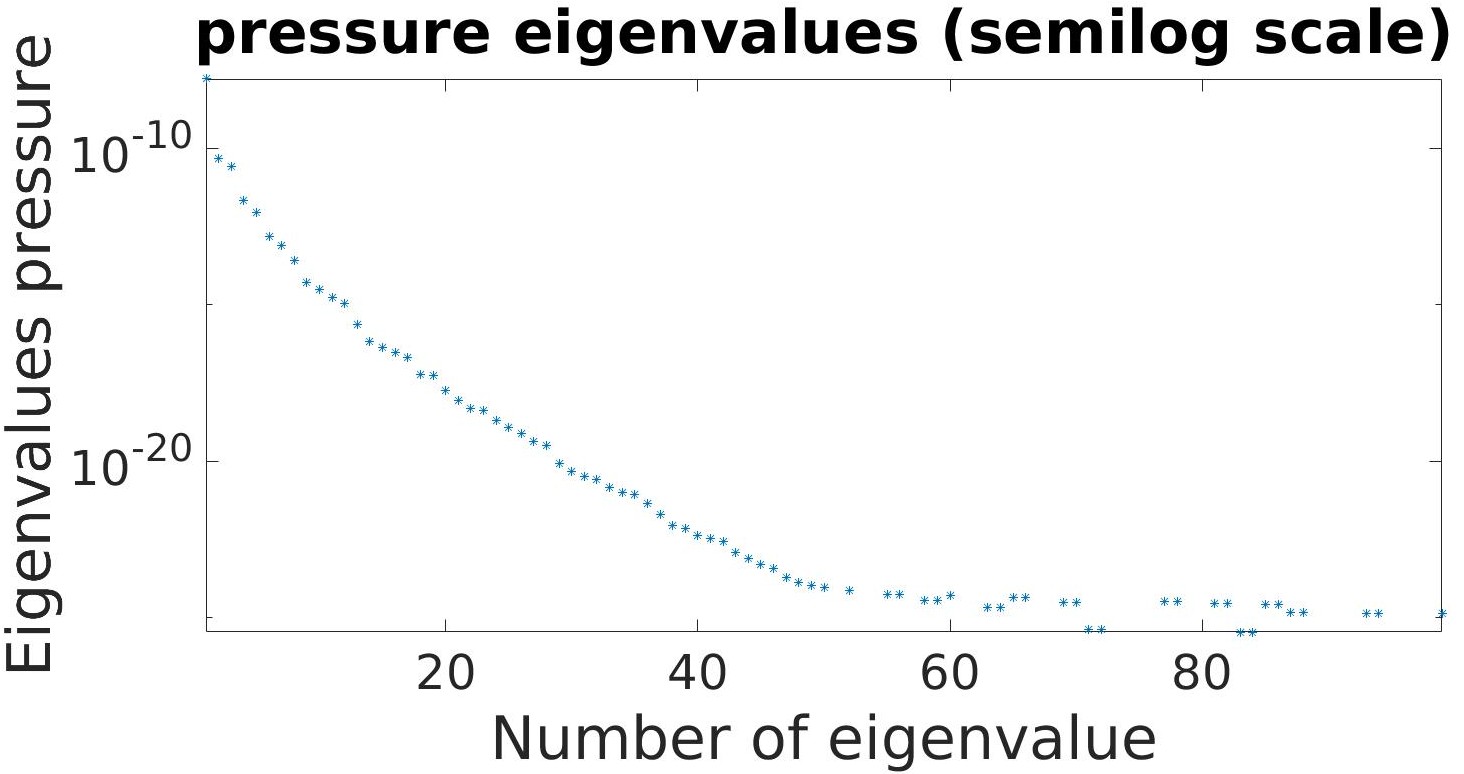}
\caption{Pressure eigenvalues (semilog scale)} \label{pressure_ev}
\end{subfigure}
\caption{Eigenvalue decay}\label{ev_decay}
\end{figure}

Figure \ref{dg_rb_solution_47_33} compares the solutions computed by DGM and Reduced Basis (RB) at parameter value $\mu = (0.47,0.33)$ with reduced basis of size $10$. The drop in error with respect to the increased size of the reduced basis space (Figure \ref{error_vs_basis}) is inline with the expectation based on the eigenvalue decay (Figure \ref{ev_decay}). The average speedup was $20.6$. Typically, during the offline phase, the full order system was assembled in $35.37$ seconds and was solved in $6.74$ seconds. During the online phase, the reduced system was assembled in $2.03$ seconds and was solved in $0.009$ seconds.

\section{Some concluding remarks}

As demonstrated by the numerical example, proper orthogonal decomposition can accelerate the computations involving geometrically parametrized discontinuous Galerkin interior penalty formulation while maintaining the reliability of solution above minimum acceptable limit. The paper also discussed, the specific issues related to the geometric parametrization and the affine expansion as pertaining to the discontinuous Galerkin interior penalty formulation. We expect the current work to contribute towards exploring further potentials in the field of geometric parametrization and reduced basis approach for the discontinuous Galerkin method.

\section*{Acknowledgements}

This work has been supported by the H2020 ERC Consolidator Grant 2015 AROMA-CFD project 681447 ``Advanced Reduced Order Methods with Applications in Computational Fluid Dynamics'' and COST Action: TD1307.

\newpage
\bibliographystyle{spbasic}
\bibliography{references}

\end{document}